\documentclass[times,11pt]{article}
\usepackage{amssymb} 
\usepackage{amsmath} 
\usepackage{latexsym} 
\usepackage{epsf} 
\usepackage[all]{xy}  
\usepackage{lscape} 
\usepackage{verbatim}    
\usepackage{times}  
\usepackage{color}
 
\newtheorem{Th}{Theorem}[section]  
  
\newtheorem{proposition}[Th]{Proposition}

\newtheorem{definition}[Th]{Definition}  
\newtheorem{example}[Th]{Example}

\def\endproof{\hfill$\Box$}  
\def\qed{\ifmmode 
          $\Box$ 
         \else 
         {\unskip 
          \nobreak 
          \hfil 
          \penalty50 
          \hskip1em 
          \null 
          \nobreak 
          \hfil 
          $\Box$ 
          \parfillskip=0pt 
          \finalhyphendemerits=0 
          \endgraf} 
         \fi} 

\def\suck{\vspace{-2mm}}

\title{Ockham's razor and  reasoning about information flow}
  
\author{Mehrnoosh Sadrzadeh\footnote{Member of the research Group on Philosophy of Information  (GPI) and Information Ethics Group (IEG), University of Hertfordshire and Oxford University. Academic visitor to  Oxford University Computing Laboratory.}\\
Laboratoire Preuves Programmes et Syst\`{e}mes, \\
Universit\'{e} Paris Diderot - Paris 7 \\
{\tt mehrs@comlab.ox.ac.uk}}
\date{}  
 
\begin{document}  
\maketitle 

\begin{abstract}
What is the minimal algebraic structure to reason about information flow? Do we really need the full power of Boolean algebras with co-closure and de Morgan dual operators? How much can we weaken and still be able to reason about multi-agent scenarios in a tidy compositional way? This paper provides some answers. 
\end{abstract}

\section {Introduction}
Systems of modal logic have been applied to  disciplines of science and humanities for modeling and reasoning about concepts such as provability, time, necessity and possibility, knowledge and belief.  Each such system has its own set of axioms, specifically chosen for the domain of application it models. New application domains motivate  introduction of new axioms who may be stronger or weaker than their original peers.  Change of axioms is also motivated by development of new mathematical methods which  lead to more refined and efficient versions of the existing axioms. Introduction  of different axioms and logical systems  initiates practical and conceptual discussions on minimality issues, for instance whether or not  the set of axioms in use is the minimal such set  for the domain it promises to model, or what are the foundational structures  we want to  model in a domain and are those reflected in the  logical system we use? In this paper, we aim to further  elaborate on these issues for the case of  Epistemic Logic.  

\smallskip
\noindent
{\bf New application domains.}\\
Epistemic logics have been used by philosophers to reason about knowledge and belief, e.g. Hintikka~\cite{Hintikka} argues for the modal logic $S5$  where we have axioms that say our knowledge is truthful  and both positively and negatively introspective.   Epistemic logics have also been used by computer scientists to reason about knowledge of agents in multi-agent systems~\cite{FH, WiebeMike}. This has led to numerous different variations; a dramatic one considers logics with non-monotonic modalities to provide a solution to the problem of logical omniscience.     The need to moreover model  the interactions  among the agents and to further reason about the knowledge they  acquire  as a result of these interactions, has led to the development of \emph{logics of information flow}~\cite{FH,Ger1,Plaza,vanbenthem,vanDerhoek}.  Very roughly put, these logics are obtained by  enriching  the Epistemic logics with temporal operations.  They have been extended to  two sorted \emph{dynamic epistemic} logics which  also model dishonest interactions of cheating and lying  in~\cite{BaltagMossSolecki,BaltagMoss}. In this logic, the  propositional sort is an $S5$ Epistemic Logic with a new 'possibly wrong' belief modality, which is only  conjunction preserving, and  the dynamic sort is a  PDL-style Dynamic logic, which has linear operations of sequential composition and choice on actions and the modality  induced by it  is not  a usual closure type operator.  

\smallskip
\noindent
{\bf New mathematical models.}\\
Algebraic methods have been used by mathematical logicians to prove meta-theorems such as decidability and completeness for modal logics. For instance,  McKinsey in 1941 and  Jonsson and Tarski in 1951~\cite{Jonsson1,Jonsson2} used Boolean Algebras with operators to prove decidability and completeness  of $S2$ and $S4$. Through the work of Eilenberg and Mac Lane, algebraic systems have been generalized to  categorical structures whose operational and compositional nature lend themselves to  easier  proof theoretic implementations of logical systems. The compositionality of the categorical approach provides  more refined ways of defining operators on the logical systems. For example, the usual co-closure modalities of modal logics  $S4$ and $S5$ can seen as a  decomposition of a pair  of adjoint maps $(f \dashv g)$ whose composition $g \circ f$ will provide us with an operation which is a co-closure. The $f$ and $g$ maps themselves can be seen as `weaker' modalities of the logic, in the sense that they obey less truth axioms, for example they need not in general be idempotent (axiom 4) or reflexive (axiom $T$). However, $f$ is disjunction  preserving and $g$ is conjunction  preserving and  they relate to each other via the rule of adjunction. These  equips $f$ and $g$ with a tidy  mathematical axiomatics for reasoning about more fine-grained aspect of situations, namely those in which the main modality need not be introspective and truthful.  The generality of the  categorical approach enables us to systematically weaken the base propositions of  modal logics  and for instance work in a  Heyting algebra where the  negation operator is weaker than the one in the original Boolean algebras of Jonsson  and Tarski, or simply in  a lattice where in general no negation operator is present.   We have  taken advantage of these bonuses and  have developed an algebraic/categorical semantics to  reason about information and mis-information flow~\cite{BCS,SadrThesis}. In this algebra,  the propositional logic is a complete lattice and  the epistemic and dynamic modalities  are respectively  formed from conjunction and disjunction preserving operators  that are adjoints to one another.

\medskip
The possibility of introducing different modal logics and different ways of defining modalities makes us wonder about, and thus bring into question, the  minimal set of axioms that makes each such modality and logic a necessity for the domain they try to model. Our interest lies in the application domain of reasoning about information (and mis-information) flow where one can ask:  what is a parsimonious  logic of information flow that can model interactive scenarios of  multi-agent systems?  In other words,  what is a minimal set of modal and propositional axioms  that enable us to reason about knowledge and interaction of agents in these systems? Or more profoundly, what would be  the philosophical implications of such a minimal logic?, what kind of concepts can the  logic based on these minimal axioms reason about?, and in short,  what are the \em foundational structures \em to  reason about information flow? This paper tries to provide some answers. 

We  start by applying Ockham's razor to complete Boolean modal algebras. We observe that one can define weaker modal operators on these algebras, those that are only disjunction or conjunction preserving and as a result have left and right adjoints respectively.  In the presence of a Boolean negation, another nice  connection shows up: the De Morgan duals of these weaker modalities also become adjoints to each other. Experience shows that  most applications of Epistemic logic do not need all four of these modalities  and only  use two of them. The traditional approaches rely on negation and work on a De Morgan dual pair of these modalities. We propose a  new thesis and propose to work with an adjoint pair instead. Thus our base propositional algebra need not have a Boolean negation: it can be a Heyting Algebra or it can have no negation at all and just be a complete lattice. Interestingly enough,  fixed points can be defined for our adjoint pair of modalities  by closing them under composition and disjunction and we show that the fixed points of adjoint operators are also adjoints to each other. 

On the application side,  we provide new readings for our modalities: as \emph{appearance} and \emph{information} of agents about the reality. These are weaker than the usual \emph{knowledge} and \emph{belief} interpretations of Epistemic logics. But, and as we shall demonstrate, we can ask for extra conditions on them to re-gain the traditional   modalities  of systems $K, S4$ and $S5$. However, it would  not be very straightforward to obtain our  modalities from their stronger peers. We  then move towards the dynamic applications and extend our minimal logic with dynamic modalities  and  show, by means of examples,  how  swiftly we can prove more fine-gained and  more interesting epistemic properties of multi-agent scenarios by means of unfolding the adjunctions. Proving these properties  enables us to reason about how the \emph{information} of agents changes (not necessarily truthfully) as a result of their communication and based on their \emph{appearances}.  To model the interactions,  we first add a pair of adjoint modalities to model the temporal \emph{previous} and \emph{next} states of the system. Then, in a second incremental step,  we index   these modalities with labels. The labels stand for  actions of  multi-agent scenarios and enable us to model  what specific actions evolved the system into its  next state.  Finally, we observe that our index set, that is the  set of interactions, is more than just a plain set and admits both a monoid and a sup-lattice structure, In short, it can be seen as a \emph{quantale} with composition and non-deterministic choice of actions. Based on this observation, we end by proving how  \emph{epistemic systems} of~\cite{BCS,SadrThesis}  are obtained from  our incrementally developed \emph{real action epistemic algebras} by restricting their agents to the  \emph{optimistically paranoid} ones.

This paper can also be seen as  a deductive take on the algebraic semantics of dynamic epistemic logic as presented in~\cite{BCS,SadrThesis}.   We  demonstrate how the full structure is put together operation  by operation, and what new aspects of application  are modeled by each operation. All along, we follow  the same parsimonious strategy for  both epistemic and action modalities, our strategy  shows that the reliance of  traditional modal and epistemic logics  on negation (classical and intuitionistic) can be waived  by  using adjoint operators instead. The theoretical study of this minimal modal algebra, its free construction and equational theory constitutes future work. 

\section{Ockham's razor and reasoning about information}
Intuitively,  a  Boolean algebra can be seen  as a propositional logic in the following way
\begin{itemize}
\item Elements of the algebra $b_1, b_1 \in \cal B$ are logical propositions, 
\item their join $b_1 \vee b_2$ is the logical disjunction, 
\item their meet $b_1 \wedge b_2$ is the logical conjunction, 
\item and the partial order between them $b_1 \leq b_2$ is the logical entailment.  
\end{itemize}
The definition of  a complete Boolean algebra~\cite{Priestely} is as follows
\begin{definition}\em
A \emph{complete Boolean algebra} ${ \cal B} = (B, \bigvee, \neg)$ is  a distributive  complete lattice $(B, \bigvee)$ with   a negation operation,  defined by the following axioms
\[
b \wedge \neg b = \bot, \qquad b \vee \neg b = \top
\]
\end{definition}
A complete Boolean algebra has all joins  $\bigvee_i b_i$, in particular the empty one $\bigvee \emptyset = \bot$, as a result it also has all meets $\bigwedge_i b_i$,  in particular the empty one $\bigwedge \emptyset = \top$.  The negation is involutive, that is $\neg \neg b = b$.

A complete Boolean algebra\footnote{For simplicity  of presentation we work with complete lattices, so that for every join preserving operator there exists a right adjoint. An alternative would be to put aside the completeness criteria and instead ask for existence  of adjoints for each join (meet)  preserving operator.}  is endowed with operators  that satisfy  certain properties,  to obtain an algebraic [classical] modal logic. In this setting,   unary operators will stand for modalities of the logic. We define a \emph{classical modal algebra}  as follows 
\begin{definition}\em
A \emph{classical modal algebra} ${\cal B} = (B, \bigvee, \neg, f) $ is a complete Boolean algebra $(B, \bigvee, \neg)$ endowed with a join preserving  operator  $f \colon B \to B$, that is 
\[
f(\bigvee_i b_i) = \bigvee_i f(b_i), \qquad \text{in particular} \ f(\bot) = \bot
\]
\end{definition}
So far we have one modality, that is the  $f$ operator, which preserves the disjunctions of the logic. But recall that in every classical  modal algebra $(B, \bigvee, \neg, f)$,  the join preserving operator $f\colon B \to B$ has a  de Morgan dual $g \colon B \to B$ defined as
\[
g(b) := \neg f(\neg b)
\]
satisfying
\[
g(\bigwedge_i b_i) = \bigwedge_i g(b_i), \qquad  \text{in particular} \  g(\top) = \top
\]
So we obtain another modality, that is the $g$ operator, which preserves the conjunctions of the logic. The categorical methods remind us  that in a classical modal algebra $(B, \bigvee, \neg, f)$, a join preserving operator $f \colon B \to B$ has a  meet preserving Galois right adjoint~\cite{JoyalT}, denoted by $f \dashv f^*$ defined as
\[
f^*(b) := \bigvee \Big\{b' \in B\mid f(b') \leq b \Big\}
\]
Moreover and in a similar way,  the  de Morgan dual of $f$, abbreviated as $g$,  has a  join preserving Galois left adjoint,  denoted by $ g^* \dashv g$ and defined as
\[
g^*(b) := \bigwedge \Big\{b' \in B\mid b  \leq g(b') \Big\}
\]
So we have obtained four modalities: two De Morgan duals and two adjoints. The adjoint operators $f \dashv f^*$ and $g^* \dashv g$ satisfy the following rules
\[
f(b) \leq b' \quad \text{iff} \quad  b \leq f^*(b'), \qquad g^*(b) \leq b' \quad \text{iff} \quad  b \leq g(b')
\]
As  a consequences, the following hold for the composition of adjoints
\[
f (f^*(b)) \leq b, \quad b \leq f^*(f(b)), \quad g^*(g(b)) \leq b, \quad b \leq g(g^*(b))
\]
There is an interesting cross-dependency between the pair of our De Morgan dual modalities $(f,g)$ and their adjoints $f^*$ and $g^*$, namely that the adjoints to the De Morgan duals are De Morgan duals of one another. In other words, in  a classical modal algebra the  de Morgan duality between $f$ and $g$ lifts to $f^*$ and $g^*$ as shown below
\begin{proposition}
In a classical modal algebra $(B, \bigvee, \neg, f)$ the following is true
\[
 f^*(b) = \neg g^*(\neg b)
\] 
where $g$ is the de Morgan dual of $f$, and we have $f \dashv f^*$ and $g^* \dashv g$.
\end{proposition}
{\bf Proof.} We show $f^*(b) \leq \neg g^*(\neg b)$ and  $\neg g^*(\neg b) \leq f^*(b)$. 
For the first inequality, start from the consequence of adjunction $f(f^*(b)) \leq b$,   by anti-tonicity of negation it follows that  $\neg b \leq \neg f(f^*(b))$, by involution of negation this is equivalent to  $\neg b \leq \neg f(\neg \neg f^*(b))$, by de Morgan duality between $f$ and $g$ this is equivalent to  $\neg b \leq g(\neg f^*(b))$, by the adjunction rule between $g$ and $g^*$ this is iff  $g^*(\neg b) \leq \neg f^*(b)$, which implies $f^*(b) \leq \neg g^*(\neg b)$ by anti-tonicity of negation. Proof of the other inequality is similar. 
 \endproof

Thus in a  complete Boolean algebra ${\cal B} = (B, \bigvee, \neg)$ asking for $f \colon B \to B$ immediately provides us with 3 other  maps $f^*, g, g^*$, which form two pairs of adjoint operators  and two pairs of de Morgan dual operators:
\[
\Big(f \dashv f^*, \quad g^* \dashv g\Big ), \qquad  \Big(g(-) := \neg f(\neg -), \quad g^*(-) := \neg f^*(\neg -)\Big)
\]
If we weaken the base algebra from a Boolean algebra $BA$ to a Heyting algebra $HA$ and thus obtain an intuitionistic modal algebra, a join preserving $f$ operator gives rise to three other operators, in the same way as in a Boolean algebra. However,  because  the Intuitionistic negation  defined as $\neg l := l \to \bot$, for $\to$ the left adjoint to $\wedge$, is not involutive,  $f^*$ and $g^*$ will not  be de Morgan duals any more.  Thus by weakening the base algebra we also obtain  weaker connections between the  operators on the base. If we continue this weakening and reduce  the base algebra  to a distributive complete lattice $DL$, there is no negation operator present and the $f$ map only gives rise to an $f^*$, the same is true in a complete lattice $L$.  The relation between the base algebra and the operators defined on it is depicted in the table below

\bigskip
\begin{center}
\begin{tabular}{|c|c|c|}
\hline
{\bf Negation} & {\bf De Morgan dual modalities} & {\bf Adjoint modalities}\\
\hline
\hline
{\tt Classical negation}&$(BA, f, g)$ & $(BA, f \dashv f^*)$\\
\hline

&$(BA, f^*, g^*)$  & $(BA,  g^* \dashv g)$\\
\hline

&&\\
\hline

{\tt Intuitionistic  negation} & $(HA, f, g)$ & $(HA, f \dashv f^*)$\\
\hline

& -& $(HA, g^* \dashv g)$\\
 \hline
 
& & \\
 \hline
 
 {\tt No negation} & - & $(DL, f \dashv f^*)$\\
 \hline
 
 &-&- \\
 \hline
 
{\tt No negation} &- &  $(L, f \dashv f^*)$ \\
\hline
&-&-\\
 \hline
\end{tabular}
\end{center}

\bigskip
The modal logic based on the algebra of the  last line of the table is weaker than the modal logics based on the algebras of its above lines, in the sense that  it  asks for  the least set of axioms from its base algebra and operators. In this case, the base algebra is a complete lattice with only one join preserving operator. We  focus on   this algebra  as our minimal modal algebra. More formally, we have
\begin{definition}\em
An \emph{adjoint modal algebra} denoted by $(L, f \dashv f^*)$ is 
a complete lattice $L$ endowed with a join preserving map $f \colon L \to L$.
\end{definition}
The operators $f$ and $f^*$ and the adjunction $f \dashv f^*$ between them can be used to define other pairs of adjoint maps on the base algebra. For example,  closing both $f$ and $f^*$  under composition  and disjunction provides us with a pair of interesting operators. These closed maps  can be seen as  special \emph{fixed point} operators, which   will stay adjoint to each other, via the following result
\begin{proposition}\label{fixedpoint}
In any adjoint modal algebra $(L, f \dashv f^*)$, the following are true
\begin{itemize}
  \item $f^i \dashv f^{*\, i}\,, \quad$   $\forall i \in \mathbb{N}$, where $ f^i = \underbrace{f \cdots f}_{i}$ stands for  $i$ times self composition of  $f$.
  \item $\bigvee_{i=1} f^{i} \dashv \bigwedge_{i=1} f^{*\, i}$
\end{itemize}
 \end{proposition}
 
 \medskip
 \noindent
{\bf Proof.} $f^i \dashv f^{*\, i}$ is equivalent to  $f^i(l)  \leq l' \  \text{iff} \ l \leq f^{*\, i}(l')$, which follows by  $i$ times applying $f (l) \leq l' \ \text{iff} \  l \leq f^*(l')$.  Similarly,  $\bigvee_{i=1} f^{i} \dashv \bigwedge_{i=1} f^{*\, i}$ is equivalent to $\left (\bigvee_{i=1} f^{i} \right ) (l) \leq l'\ \text{iff}\  l \leq  \left (\bigwedge_{i=1} f^{*\, i}\right )(l')$, which follows  from the definitions of arbitrary meets and joins applied to item \hbox{one.\endproof}

\medskip
Operators of the first item above are closed under composition and can be seen as a pair of adjoint fixed point operators.  Operators of the second item above are moreover closed under disjunction and conjunction respectively and can be seen as  a pair of adjoint least and greatest fixed points. One can make these reflexive by starting the range of  $i$ from 0.  

We make our modal  algebra more suited  for epistemic applications by  considering a family of join preserving operators, instead of just one, and thus obtain a multi-modal algebra, defined below
\begin{definition} \em
A \emph{multi-agent adjoint modal algebra} (MAMA) denoted by $(L, f_A \dashv f^*_A)_{\cal A}$ is a complete  lattice $L$ endowed with a family of join preserving maps $\{f_A\}_{A \in {\cal A}}\colon L \to L$.
\end{definition}

We refer to this algebra as an \emph{epistemic algebra} and provide epistemic interpretations for its modalities:  $f_A(l)$ is interpreted as  `appearance of  proposition $l$ to agent $A$'. That is  
\begin{center}
\fbox{$f_A(l)$ is all the propositions that \emph{appear} to agent $A$ as possible or true when in reality $l$ is true.}
\end{center}
Here are some explanatory examples of our notion of \emph{appearance}:
\begin{itemize}
\item If $f_A(l)=  l$ then the appearance of agent $A$  about reality is the reality itself, so $A$'s appearance is totally compatible with reality. 
\item If $f_A(l)=  \top$ then all the propositions of the logic appear as possible to agent $A$, in other words, he has no clue about what is going on in reality.
\item If $l \leq f_A(l)$, for instance when  $f_A(l) = l \vee l'$, then reality appears as possible to agent $A$, although he cannot be sure about it, since $l'$ also appears equally possible to him.
\end{itemize}
In each case above, we can also  talk about \emph{information} of agent $A$, in the following lines
\begin{itemize}
\item If $f_A(l)=  l$ then $A$'s information about reality is the reality itself, so $A$ is well informed or has truthful information.
\item If $f_A(l)=  \top$ then $A$ has no information at all about reality.
\item If   $f_A(l) = l \vee l'$, then $A$'s information  about reality  includes the reality, but is weaker than it.
\end{itemize}
Based on the above intuitions,  we  use the  left adjoint  $f^*_A(l)$ to define our notion of \emph{information} and  read   it  as   `information of agent $A$ about proposition $l$', or more propositionally as follows
\begin{center}
\fbox{ $f^*_A(l)$ is read as  `agent $A$ is informed that proposition $l$ holds'.}
\end{center}
Now we can apply the adjunction rule
 \[
f_A(l) \leq l  \quad \text{iff} \quad l \leq f^*_A(l')
\] 
to produce equivalent information formulae for each appearance case above: 
\begin{itemize}
\item If $f_A(l)=  l$ then $l \leq f^*_A(l)$, so   $l$ implies that $A$ has truthful information about $l$.
\item If $f_A(l)=  \top$ then $l \leq f^*_A (\top)$, so $l$ implies that $A$ has no information about $l$.
\item If   $f_A(l) = l \vee l'$, then $l \leq f^*_A(l \vee l')$, so $l$ implies that  $A$ is informed that either $l$ or another proposition $l'$  hold in reality.
\end{itemize}
For more examples on appearances consider the following comparisons
\begin{itemize}
\item If $f_A(l) \leq f_B(l)$ then agent $B$ is more uncertain about $l$ than agent $A$, since  more propositions appear as possible to him. So we can say that $A$ is  more informed or has more information about proposition $l$ than agent $B$. An example would be  when $f_B(l) = l \vee l'$ where as $f_A(l) = l$, clearly $l \leq l \vee l'$ and  $l \vee l'$ stands for two possibilities for agent $B$ as opposed to the only one possibility, that is $l$,  for agent $A$.
\item If $f_A(l) \leq f_A(l')$, then agent $A$ is more uncertain about $l'$ than about $l$, thus he is more informed about $l$ than about $l'$.  An example would be when $f_A(l) = l$ and $f_A(l') = l \vee l'$, so whenever  $l$ is true in reality, $A$ is informed that this is the case, but when $l'$ is true in reality,  his information does not tell anything useful to him, since he  cannot distinguish between $l$ and $l'$, both appear to him as equivalently possible.
\end{itemize}

 The above notions of \emph{appearance} and \emph{information} are based on weaker modalities  than those of the usual Epistemic logics. They  provide new readings for the modalities; readings that stand for new concepts that Epistemic logics did not account for before. However, in our weaker system, we can define the stronger epistemic notions of   other logics. For instance, the knowledge modality of system $K$ can now be  described as `truthful information' and defined by  
  \[
  K_A \, (l) := f_A^*(l) \wedge l
  \]
  So we have
  \begin{center}
  \fbox{$K_A\, (l)$  is read as  '$A$ has truthful information that $l$'.}
  \end{center}
If  the appearance maps are  weakly idempotent and decreasing (i.e.~  weak co-closures), then one obtains the  knowledge of system $S4$.
\begin{proposition}\label{S4}
In a MAMA $(L, f_A \dashv f^*_A)_{\cal A}$, if  we have $f_A(l) \leq l$ and $f_A f_A(l) \leq f_A(l)$  then the following hold
\begin{itemize}
\item  $f_A^* (l) \leq f_A^* f_A^* (l)$ 
\item  $K_A(l)  \leq K_A K_A(l)$ \ and \ $ K_A(l) = l$,  \ for  \ $K_A(l) := f_A^*(l) \wedge l$
\end{itemize}
\end{proposition}
{\bf Proof.}  For the first one, by the corollary of adjunction we have $f_A f_A^* (l) \leq l$, from this by weak idempotence of $f_A$ and transitivity we have $f_A f_A f_A^* (l) \leq f_A f_A^*(l) \leq l$ and thus it follows that $f_A f_A f_A^* (l) \leq l$, which by adjunction is equivalent to $f_A^* (l) \leq f_A^* f_A^* (l)$. For the second one, we have to show $f^*_A(l) \wedge l \leq f^*_A(f^*_A(l) \wedge l) \wedge (f^*_A(l) \wedge l)$, which is equivalent to $f^*_A(l) \wedge l \leq f^*_A(f^*_A(l) \wedge l)$, that is $f^*_A(l) \wedge l \leq  f_A^* f_A^*(l) \wedge f_A^*(l)$,  which follows from item one and the adjunction equivalence of decreasing property of $f_A$, that is $l \leq f_A^*(l)$. The third one easily follows from $f_A(l) \leq l$, which is  by adjunction equivalent to $l \leq f^*_A(l)$ and by definition of meet we obtain $l \wedge f^*_A(l) = l$, which is nothing but $K_A(l) = l$. 
\endproof

\medskip
When   $L$  is a complete Boolean algebra,  belief  is defined as  the  de Morgan dual of $K_A$, that is $B_A(l) := \neg K_A(\neg l)$.  In this setting, knowledge  of the system $S5$ is obtained by asking for the weak idempotence and decreasing of appearance maps. 
\begin{proposition} In a MAMA $(L, f_A \dashv f^*_A)_{\cal A}$, if  $L$ is a complete Boolean algebra  $(L, \bigvee, \neg)$ and we have  $f_A(l) \leq l$ and $f_A f_A(l) \leq f_A(l)$  then it follows that 
\begin{itemize}
\item $\neg K_A \leq K_A( \neg K_A(l))$
\item  $ K_A(l) = l$.
\end{itemize}
\end{proposition}
{\bf Proof.} Similar to the proof of proposition~\ref{S4}.

\medskip
\noindent
Similar to proposition~\ref{fixedpoint}, we define adjoint fixed points for our indexed modalities as below
\begin{proposition}\label{multifixedpoint}
In any  MAMA  $(L, f_A \dashv f^*_A)_{\cal A}$ \ for \ $\beta \subseteq {\cal A}$ the following are true
\begin{itemize}
\item $f_{\beta}  \dashv  f^*_{\beta}\,,\quad$ for     $f_{\beta}  :=  \bigvee_{B \in \beta} f_B$ \ and \ $f_{\beta}^* :=  \bigwedge_{B \in \beta} f_B^*$.
  \item $f^i_{\beta} \ \dashv \  f^{*\,i}_{\beta}$
  \item  $ \bigvee_{i = 1} f^i_{\beta} \ \dashv \  \bigwedge_{i = 1} f^{*\,i}_{\beta}$\\
\end{itemize}
 \end{proposition}
  {\bf Proof.} For the first direction of the first one assume $f_{\beta} (l) \leq l'$,  by definition of join it follows that $f_B(l) \leq l'$  for all $B \in \beta$, by adjunction  this is iff $l \leq f^*_B(l')$ for all $B \in \beta$,  by definition of meet it follows that $l \leq f^*_{\beta} (l')$.  Proof of the other direction is similar.   The second one follows from $i$ times unfolding the first one, and the third one from the first two. \endproof

    \medskip
    \noindent
 These group maps have sensible interpretations in an epistemic context, for example $f_\beta (l)$ can be  read as  'the appearance of $l$ to all the agents in group $\beta$', similarly $f^*_{\beta} (l)$ can be read as   'the  shared information of agents in $\beta$ about $l$', or more propositionally as
 \begin{center}
 \fbox{$f^*_{\beta} (l)$ is read as 'all the agents in $\beta$ are informed that $l$ holds.}
 \end{center}
  The former contains  the collection or the union of appearances of agents in $\beta$ about the same proposition $l$,  and can be read as `accumulated appearance'. The latter contains the common part or the intersection of information of agents in $\beta$ about the same proposition $l$ and can be  read as `shared information'. 
Closing the \emph{shared information} under composition and conjunction (the third item in proposition~\ref{multifixedpoint}) provides us with the notion of \emph{common information}, which is the infinite nested information of agents about one another's information:
\begin{center}\fbox{\begin{minipage}{14cm}\begin{center}
$\bigwedge_{i = 1} f^{*\,i}_{\beta} (l)$  is read as 
`all the agents in group $\beta$ are informed that $l$, and are also informed that everyone in the group is informed that $l$, and so on $\cdots$'. \end{center}\end{minipage}}
\end{center}
The notion of `common knowledge among the group $\beta$' in system $K$  is obtained by starting the index $i$ from 0 rather than 1, that is 
\[
CK_\beta := \bigwedge_{i = 0} f^{*\,i}_{\beta} = l \wedge \bigwedge_{i = 1} f^{*\,i}_{\beta}
\]
In other words, common knowledge among agents in  group $\beta$ can be defined in terms of their common  information as follows
\begin{center}
\fbox{Agents in $\beta$ have common knowledge that $l$ iff they have truthful common information that $l$.}
\end{center}

\medskip
Although, weaker than the knowledge and belief modalities of Epistemic logics, appearance and information modalities can also be used to model  epistemic applications and to prove weaker properties about them.  However, and as we will see in the next section, this weakness becomes a necessity while reasoning about mis-information.   
\begin{example} 
Consider the following coin toss  scenario:
in front of agents $A$ and $B$, agent $C$ throws a coin and covers it in his palm. We consider  a  MAMA containing  propositions $H, T \in L$. Appearances are set according to uncertainty of agents
\[
f_A(H) = f_A(T) = H \vee T
\]
We show that the information $A, B$ and $C$  have is that the coin   is either heads or tails, for example
\begin{eqnarray*}
H &\leq& f^*_A (H \vee T)\\
H &\leq&  f^*_A\,   f^*_B (H \vee T)\\
H &\leq&  f^*_B\,   f^*_A\, f^*_B (H \vee T)
\end{eqnarray*}
Consider the second property,  by the adjunction rule  it holds iff we have  $f_A(H) \leq  f^*_B (H \vee T)$, 
by assumptions on  $f_A$ this is  equivalent to $H \vee T \leq f^*_B (H \vee T)$. By the adjunction rule this holds iff we have $f_B(H \vee T ) \leq H \vee T$, now 
since $f_B$ is join preserving this is equivalent to  $f_B(H) \vee f_B(T) \leq H \vee T$, which is, 
by  assumptions on $f_B$, equivalent to  $ (H \vee T) \vee (H \vee T) \leq H \vee T$, which holds by the definition of $\vee$.  The proofs of other cases are similar.\end{example}

\section{Ockham's razor and reasoning about flow of information}
To reason about flow of information, we  add another modality to our epistemic algebra: the  action modality. This  will enable us to  prove more properties about scenarios: before we were able to prove   that   agents have some information, now we can show how their acquired this information, that is  how  their initial information  got updated as a result of  some communication action  taking place among them. Epistemic algebras could only reason about information, the question is how to enrich them in a minimal way such that they can also reason about communication.
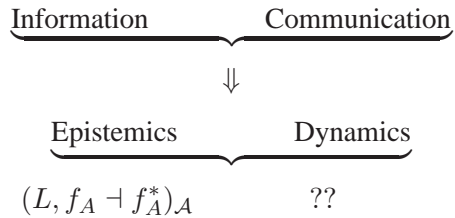
\begin{center}
\begin{minipage}{15cm}
\[
\underbrace{\text{Information} \qquad  \qquad \text{Communication}}
\]
\[
\Downarrow
\]
\[
\underbrace{\text{Epistemics} \qquad  \qquad \text{Dynamics}}
\]
\[
\hspace{-1.4cm}(L, f_A  \dashv f^*_A)_{\cal A} \qquad \qquad ??
\]
\end{minipage}
\end{center}

\medskip
\noindent
The resulting logic  is  obtained by  endowing our MAMA with a new operator to stand for dynamics. The reasoning power of this algebra is increased by asking  the new operator to weakly permute with the existing epistemic operators. The definition of the new algebra is as follows
\begin{definition}\em
A  \emph{temporal epistemic algebra}  denoted by $(L, f_A \dashv f^*_A, h \dashv h^*)_{{\cal A}}$
 is a multi-agent adjoint modal algebra $(L, f_A \dashv f^*_A)_{{\cal A}}$ endowed with a join preserving map $h \colon L \to L$, such that the following permutation holds
\[
 \quad f_A\,  h (l)  \leq h\,  f_A(l)\\
\]
\end{definition}
We read $h^*(l)$ as 
\begin{center}
\fbox{`In the next state of the system $l$ holds'.}
\end{center}
Hence $h^* f^*_A (l)$ is read as 
\begin{center}
\fbox{`In the next state of the system agent $A$ gets informed that $l$ holds.}
\end{center}
 Similarly,  $h(l)$ is read as `in the previous state of the system $l$ held'.
The permutation axiom of the algebra is a weak permutation between the two operators and  demonstrates a preservation  or no-miracle condition on the information: if an agent obtains some information in the next state of the system, it should be the  case that this information existed in the  system  previously, thus  the acquired information should somehow be implied by the previous information. In other words, information is not generated and cannot be destroyed freely and without a cost, it can only be accumulated.   This axiom is similar  to the \emph{appearance-update} axiom of~\cite{BCS,SadrThesis} and  also corresponds to a weaker  version of the \emph{action-knowledge} axiom of~\cite{BaltagMoss, BaltagMossSolecki}.  A similar axiom can also be found in the Epistemic Temporal Logic of~\cite{FH} in the name of \emph{perfect recall}. 

According to proposition~\ref{multifixedpoint}, a temporal fixed point operator  can be defined as  $\bigwedge_i h^{*\, i} (l)$   and interpreted as follows
\begin{center}
\fbox{'Eventually in some future state of the system proposition $l$ holds'.}
\end{center}

Rather than reasoning about whether  an information property holds  in the next state of the system, it would be more sensible to  name and reason about  the action that led the system to its next state, the action that caused the information property to hold in the next state of the system.    To do so,   we endow our  temporal epistemic algebra with a family  of operators that are indexed over a set of actions.  The passage from one temporal operator to a family of action operators is similar to the passage from the mono-modal epistemic algebras to the multi-agent ones. The new setting is defined as follows
\begin{definition}\em
An   \emph{action epistemic algebra}   denoted by  \ $(L, f_A \dashv f^*_A, h_a \dashv h^*_a)_{{\cal A}, Act}$ \  is a multi-agent adjoint modal algebra  \ $(L, f_A \dashv f^*_A)_{{\cal A}}$ \ endowed with a family of join preserving maps \hbox{$\{h_a\}_{a \in Act} \colon L \to L$}, such that 
\[
f_A\,  h_a (l)  \leq h_a\,  f_A(l)
\]
\end{definition}
We  read $h^*_a(l)$ as 
\begin{center}
\fbox{`After action $a$ proposition $l$ holds'.}
\end{center}
Similarly,   $h^*_a\,   f^*_A (l)$ is read as
\begin{center}
\fbox{`After action $a$ agent $A$ gets informed that  $l$ holds'.} \end{center}
 The fixed point of the action operator $\bigwedge_i h_{\alpha}^{*\,i}$, for $\alpha \subseteq Act$ is interpreted as 
 \begin{center}
 \fbox{`Eventually after the actions in $\alpha$ proposition $l$ holds'.}
 \end{center}
 We end this section by describing  two  restrictions that will bring our algebras  closer to the specific application domain in mind. These restrictions  have been introduced and discussed in detail in the algebra of~\cite{BCS,SadrThesis}, and  correspond to similar restrictions in the dynamic epistemic logic of~\cite{BaltagMoss,BaltagMossSolecki}.  The main point is that  the actions that we are interested in reasoning about   are the communication actions that take place in epistemic scenarios. These do not change the facts of the world and are of the form of  announcements to a group of agents of  a propositional or epistemic content. In order to model them,  we ask for the following two axioms,   for  $a \in CAct \subseteq Act$   and $\phi \in \Phi \subseteq L$ 
 \begin{eqnarray*}
  l \in ker(a)& \quad \text{iff} \quad & h_a(l) \leq \bot\\
 l \leq \phi &\quad \text{iff} \quad & h_a(l) \leq \phi
\end{eqnarray*}
We refer to $CAct$ as the communication actions and to $\Phi$ as the 'facts' of the system. 
The first axiom says that each communication action $a \in CAct$ has a kernel $ker(a)$, which stands for its `co-content', that is all the propositions  to which the action cannot be applied.  The second axiom says that if a proposition $l$  entails a fact $\phi$, that is $l \leq \phi$ then a communication actions $a$  does not have any effect on this entailment, that is $h_a(l) \leq \phi$.

\begin{example}
As an example, consider again the coin toss scenario where agent $C$ uncovers the coin and announces: `the coin is heads'. The announcement is a communication action $a \in CAct$ that appears as it is to all the agents since it is a  public action, so $f_A(a) = f_B(a) = f_C(a) = a$. The kernel of this action is $T$, since it cannot apply when the coin has come down tails. The set of facts in this scenario is  $\{H,T\}$. 

We want to show that after this announcement the uncertainty of agents gets waived and for instance $A$ will acquire information that the coin is heads, that is 
\[
H \leq h^*_a f^*_A(H)
\]
By the adjunction rule on $h^*_a$ this   holds iff $h_a(H) \leq f^*_A(H)$, by the adjunction rule on $f^*_A$, this holds iff $f_A h_a (H) \leq H$. By the no-miracle axiom it suffices to show  $h_a f_A(H) \leq H$. By the assumptions on $f_A(H)$ this is equivalent to  $h_a (H \vee T) \leq H$. Since $h_a$ is join preserving this is equivalent to showing $h_a(H) \vee h_a(T) \leq H$. By definition of $\vee$ in a lattice it suffices to show the following two case 
\[
\begin{cases}
h_a(H) \leq H&\\
h_a(T) \leq H&
\end{cases}
\]
The first case follows since $H$ is a fact and thus $h_a(H) = H$ and in a partial order we have that   $H \leq H$. The second case follows since $T \in ker(a_H)$, which means $a_H(T) = \bot$, thus $\bot \leq H$. 
Other nested information  properties  such as the following ones are proved in a similar fashion
\[
H \leq h^*_a f^*_A f^*_C(H), \quad H \leq h^*_a f^*_A f^*_B(H)
\]
\end{example}

\section{Ockham's razor and reasoning about flow of mis-information}
In  the closer-to-real life versions of the scenarios of multi-agent systems, agents are not always honest and thus communication actions are not always truthful. We would like to be able to reason about these scenarios and model the cheating and lying  actions of dishonest agents. In order to do this,  and following the approaches of~\cite{BCS,BaltagMoss, BaltagMossSolecki,SadrThesis}, we introduce epistemic structure on actions. Similar to the epistemic structure on propositions, these will stand for  'appearances of agents about actions'.  Also similar to the epistemic structure on propositions, these are added by endowing the set of actions $Act$  with a family of appearance maps $f'_A \colon Act \to Act$, one for each agent.  The new algebras are defined below
\begin{definition}\em
A \emph{real action epistemic algebra}  denoted by  $(L, f_A \dashv f^*_A, h_a \dashv h^*_a)_{ {\cal A}, (Act, f'_A \dashv (f'_A)^*)_{\cal A}}$  is an action epistemic algebra $(L, f_A \dashv f^*_A, h_a \dashv h^*_a)_{ {\cal A}, Act}$ where the set of actions $Act$ is endowed with  a family of join preserving maps $\{f'_A\}_{A \in {\cal A}} \colon Act \to Act$, and we have
\[
f_A\,  h_a (l)  \leq h_{f'_A(a)}\,  f_A(l)
\]
\end{definition}
The no-miracle axiom now becomes a no-miracle axiom up to the appearance  of actions, that is if an agent acquires new information after an action, this information is based on the state of the system before the action and also the appearance of the agent about that action. 

\begin{example}
Consider the coin toss scenario,  we show that if $C$'s announcement was not honest and he lied about the face of the coin, that is announced heads when he saw tails, $A$ and $B$, who did not notice and neither suspect the lying,  will acquire wrong information. The lying action is a communication action $\overline{a} \in CAct$ that appears as it is to the announcer $C$, that is $f'_C(\overline{a}) = \overline{a}$, but since $A$ and $B$ do not suspect it they think it is an honest announcement that is $f'_A(\overline{a}) = f'_B(\overline{a}) = a$. The kernel of the lying action is $H$ since it could not be a lie if the coin had actually landed heads.  In this lying scenario we can show, for example, the following properties
\[
H \leq h^*_{\overline{a}}  f^*_C(T), \quad H \leq h^*_{\overline{a}}  f^*_A(H), \quad H \leq h^*_{\overline{a}}  f^*_Af^*_C(H)
\]
Properties of this and other examples, such as the muddy children puzzle with cheating and lying, are  proved using the same strategy as in the honest versions demonstrated. In the muddy children one needs to repeat the kernel argument for the number of dirty children in the puzzle minus 2. 

Consider the third property, by adjunction  on $h^*_{\overline{a}}$,  \ $ f^*_A$ \ and \ $f^*_C$ respectively, it is equivalent to $f_C\,  f_A\,  h_{\overline{a}} (H) \leq H$. By the no-miracle  axiom between $f_A$ and $h_{\overline{a}}$  it suffices to show  
\[
f_C\,  h_{f'_A(\overline{a})}\,  f_A (H) \leq H
\]
 which is equivalent to $f_C\,  h_{a}\,  f_A (H) \leq H$ since $f'_A(\overline{a}) = a$. By the no-miracle axiom this time between  $f_C$ and $h_{a}$, it suffices to show $h_{f'_C(a)} \, f_C\,  f_A(H) \leq H$, which is equivalent to $h_{a} f_C f_A(H) \leq H$ since $f'_C(a) = a$. We substitute values for $f_A$ and $f_C$ and need to show $h_{{a}}(H \vee T) \leq H$. By distributivity and definition of join this is obtained by showing two cases
\[
\begin{cases}
h_{{a}}(H ) \leq H&\\
h_{{a}}(T) \leq H
\end{cases}
\]
The second case holds since $T$ is in the kernel of ${a}$, and thus  $h_{{a}}(T)  = \bot \leq H$, the first case follows similar to the previous example and by preservation of facts. 
\end{example}

Two observations are in place here:
\begin{itemize}
\item  The family of indexed unary maps $\{h_a\}_{a \in Act} \colon L \to L$ is equivalent to the binary operation  \[
h \colon L \times Act \to L\,.
\]
\item There is some implicit structure on the set of actions:  they can be sequentially composed  $a \bullet a'$, non-deterministically chosen $a \vee a'$, and  there is a neutral action 1 in which nothing happens  $1 \bullet a = a \bullet 1 = a$. Assuming the existence of all the choices (joins) and their distributivity over the composition, permits  us to form a \emph{quantale} of actions  ${\cal Q} = (Q, \bigvee, \bullet,  1)$\footnote{This can be, for instance,  the powerset of the free monoid generated on $Act$, that is ${\cal P}(Act^*)$.}. 
\end{itemize}

The index sets of a real action epistemic algebra make the structure a bit too crowded, especially the index set of actions which is itself  indexed over the set of agents. The situation can be improved by considering instead two separate multi-agent adjoint modal algebras: one for  the propositions 
$(L, f_A \dashv f^*_A)_{\cal A}$ and another one for the actions $(Act,   f'_A \dashv (f')^*_A)_{\cal A}$ where the latter acts on the former via the binary counterpart of the $h_a$ operators, that is via the binary operation of  $h \colon L \times Act \to L$.   It is easy to show that  the equivalence mentioned in  the first observation above lifts to one between a real action epistemic algebra  and these two MAMA's. Formally speaking we have
\begin{proposition}\label{two-sorted}
A real action epistemic algebra 
\[
(L, f_A \dashv f^*_A, h_a \dashv h^*_a)_{ {\cal A}, (Act, f'_A \dashv (f'_A)^*)_{\cal A}}
\]
 is  equivalent to
  \[
\big((L, f_A \dashv f^*_A)_{\cal A}, (Act,  f'_A \dashv (f')^*_A)_{\cal A}, h \big )
\]
 whenever $h \colon L \times Act \to L$ is  the binary equivalent  of  $\{h_a\}_{a \in Act} \colon L \to L$.
\end{proposition}

\noindent
{\bf Proof.} Follows directly from the equivalence of observation 1 above. In particular the join preservation of  $h_a$, that is $h_a (\bigvee_i l_i) = \bigvee_i h_a(l_i)$ lifts to $h(\bigvee_i l_i, a) = \bigvee_i h(l_i, a)$ and the permutation between $h_a$ and $f_A$, that is $f_A h_a(l) \leq h_{f'_A(a)} f_A(l)$  lifts to $f_A h(l,a) \leq h(f_A(l), f'_A(a))$. \endproof

\medskip
So far we are able to reason about the information acquired by agents as a result of atomic actions taking place among them. The information acquired by composition of actions can also be taken care of by composing the action operators, for example $h_a h_b f^*_A (l)$ says that after doing action $a$ followed by action $b$, agent $A$ is informed that $l$.  What is missing is reasoning about the information after a choice of actions, for example to express what an agent would acquire if either action $a$ or action $b$ take place. 
If we move from the plain set of actions $Act$ to the quantale of actions $(Act, \bigvee, \bullet, 1)$, we obtain an algebraic structure on the actions which enables us  to as well   reason about  non-deterministic choices of actions.  What will happen to the appearance maps? For example, given the appearance of atomic actions  $a,b$ in $Act$, what would be the appearance of the choice of actions  $ a\vee b$, and their composition $a \bullet b$? The most natural and neutral way of extending appearance maps to choice of actions is  point-wisely, that is making  the appearance of the non-deterministic choice be equal to  the choice of the appearances 
\[
f'_A(\bigvee_i a_i) = \bigvee_i f'_A(a_i) 
\]
 How about with regard to the sequential composition and its unit? This depends on what kind of agents do  we want to model.  For instance,  we may decide to consider it possible for our agents to be \emph{paranoid}, that is,  when  nothing is happening in reality, it appears to them that  something is happening. In this case, we do not  need to ask for any extra inequalities between 1 and $f'_A(1)$.   Since, for example, $f'_A(1)$ can be equal to any action $a$, and in general there is no  order relation between 1 and an arbitrary  action $a$. However,  this may be a bit too strong of an assumption, we can weaken it by asking the agents to be \emph{optimistically paranoid}. That is,  when nothing is happening in reality, it appears to them that   either nothing is happening or something is happening. In other words,  we ask that  appearance to all agents of the action in which nothing happens always include it, that is
\[
1 \leq f'_A(1)
\]
 It is then easy to show (see~\cite{BCS,SadrThesis}) that this inequality will lead us to  an inequality on the appearance of a sequential composition, that is $f'_A (a \bullet b) \leq f'_A(a) \bullet f'_A(b)$. 
 There is a third possibility and that is when the agents  are not paranoid at all.  In other words,  whenever nothing is happening in reality,  it appears to them  that nothing is happening. So their  appearance of 1 is  equal to 1
 \[
 1 = f'_A(1)
 \]
 This inequality will force the appearance of the sequential composition to be equal to the sequential composition of the appearances, that is. It will also force our permutation axiom to be equality rather than inequality, that is 
 \[
 f'_A (a \bullet b) = f'_A(a) \bullet f'_A(b)\,, \qquad f_A h_a(l) = h_{f'_A(a)} f_A(l)
 \]
However, since the goal of this paper is to stay minimal in the axioms of the algebra and that inequality is weaker than equality, it is reasonable to  work with  the inequality versions of axioms  and assume that our agents are \emph{optimistically paranoid}.   A more detailed discussion of these and other attitudes of agents and their relation to  axioms of the algebra is well in place, but out of the limits of the current paper. 
 
Let us end by defining the notion of a  quantale endowed with appearance maps for optimistically paranoid agents and show how it will help us relate our system to the other algebra of  information and mis-information  flow.
\begin{definition}\em
An \emph{epistemic quantale} denoted by $(Q,  f'_A \dashv f'^*_A)_{\cal A}$ is a multi-agent adjoint modal algebra  where $Q$ is a quantale and moreover we have
\[
  1 \leq f'_A(1)\,, \qquad f'_A (a \bullet b) \leq f'_A(a) \bullet f'_A(b)
\]
\end{definition}

\noindent
The  \emph{epistemic systems} of~\cite{BCS,SadrThesis} are  obtained from the  two sorted structure of proposition~\ref{two-sorted}  as follows 
\begin{proposition}\label{episquant}
The pair  $\big ( (L, f_A \dashv f^*_A)_{\cal A}, (Q,  f'_A \dashv f'^*_A)_{\cal A}, h\big)$
is an \emph{epistemic system} whenever    $(L, f_A \dashv f^*_A)_{\cal A}$ is a MAMA, $(Q,  f'_A \dashv f'^*_A)_{\cal A}$ is an epistemic quantale, the pair is equivalent to a real action epistemic algebra $(L, f_A \dashv f^*_A, h_a \dashv h^*_a)_{ {\cal A}, (Act, f'_A \dashv (f'_A)^*)_{\cal A}}$ and moreover $h$ satisfies the following
\[
h(l, \bigvee_i a_i) = \bigvee_i h(l, a_i), \quad h(l, 1) = l, \quad  h(l, a\bullet b) = h(h(l, a), b)
\]
\end{proposition}
{\bf Proof.} Follows from definition~\ref{episquant} and    proposition~\ref{two-sorted}. \endproof

\section{Conclusion}
We have presented a minimal  algebraic modal logic  where modalities  are not necessarily  positively or negatively introspective, that is they do not in general obey axioms 4 and 5 of modal logic,  neither are  they in general  truthful, that is  obey axiom $T$. The only condition on them is preservation of  disjunctions or conjunctions of their base propositional setting. The propositional setting is also weak: it  has neither implication nor negation, and is not necessarily distributive.  Lack of negation means that our modalities are not de Morgan duals, but they are connected to each other in a weaker sense and as adjoints. These minimal modalities can be interpreted as new modes such as `information' and `appearance', from which belief and  truthful knowledge can be derived.  We have defined fixed point operators for  these modalities such that the pair of fixed points are also adjoints. The applicability  of our logic is demonstrated via examples of epistemic scenarios.  This weak setting can be extended to also model the flow of information, be it caused by  the passage of time or by application of  actions. Actions can have some extra structure on them to model cheating and lying, sequential composition and non-deterministic choice. All of these can be modularly  added to the weak modal algebra we started with. At the end, we show how restricting our agents to the \emph{optimistically paranoid}  ones allows us to obtain the structure of an \emph{epistemic system},  developed in previous work as the algebraic semantics of  Dynamic Epistemic Logic.

One needs to study the universal algebraic properties of our weak modal algebras in the lines of~\cite{Yde}; that  if  they have an equational theory, how can they be freely generated, what does their relational semantics look like, how to develop a Stone-like duality for them, etc. One possible challenge  might lie in  the rule of adjunction $f(l) \leq l'$ {iff} $ l \leq f^*(l')$, 
which  is not an equation. The equations are obtained from composing the adjoints,  for example $f \circ f^*(l)  \leq l$ and $ (f \circ f^*)^2 = f \circ f^*$, but those are not of rank 1. 

\end{document}